\documentclass[final]{siamltex}
\usepackage{amsmath, amsfonts, amssymb}
\usepackage{latexsym}

\newfont{\bb}{msbm10}

\newtheorem{thm}{Theorem}[section]
 \newtheorem{cor}[thm]{Corollary}
 \newtheorem{lem}[thm]{Lemma}
 \newtheorem{prop}[thm]{Proposition}
 
 \newtheorem{con}[thm]{Conjecture}

\newcommand{\iy}{\infty}
\newcommand{\bC}{{\bf C}}
\newcommand{\bT}{{\bf T}}
\newcommand{\bZ}{{\bf Z}}
\newcommand{\bR}{{\bf R}}
\newcommand{\al}{\alpha}
\newcommand{\be}{\beta}
\newcommand{\ga}{\gamma}
\newcommand{\Ga}{\Gamma}
\newcommand{\de}{\delta}
\newcommand{\eps}{\varepsilon}

\newcommand{\om}{\omega}
\newcommand{\ph}{\varphi}

\newcommand{\si}{\sigma}
\newcommand{\tht}{\theta}
\newcommand{\n}{\|}
\newcommand{\vsk}{\vspace{0mm}}
\newcommand{\vsg}{\vspace{0mm}}

\newcommand{\Ret}{{\rm Re}\,}

\title{Norms of Toeplitz Matrices with
Fisher-Hartwig Symbols\thanks{The research of the second author was supported by
Academy of Finland Project 207048. He thanks the Department of Mathematics for its
hospitality during his visit to TU Chemnitz in spring 2006.}}

\author{Albrecht B\"ottcher\thanks{Fakult\"at f\"ur Mathematik, Technische
Universit\"at Chemnitz, D - 09107
Chemnitz, Germany ({\tt aboettch@mathematik.tu-chemnitz.de})}
\and Jani Virtanen\thanks{Department of Mathematics,
University of Helsinki,
Gustaf H\"allstr\"omin katu 2b,
FI - 00014 Helsinki,
Finland ({\tt jani.virtanen@helsinki.fi})}}

\begin{document}

\maketitle
\begin{abstract}
We describe the asymptotics of the spectral norm of finite Toeplitz matrices
generated by functions with Fisher-Hartwig singularities as the matrix dimension
goes to infinity. In the case of positive generating functions, our result
provides the asymptotics of the largest eigenvalue, which is of interest in time
series with long-range memory.
\end{abstract}



\pagestyle{myheadings}
\thispagestyle{plain}
\markboth{B\"OTTCHER AND VIRTANEN}
{TOEPLITZ MATRICES WITH FISHER-HARTWIG SYMBOLS}

\section{Introduction}\label{S1}

Let $\{a_k\}_{k \in \bZ}$ be a sequence of complex numbers and denote
by $T_n$ the $n \times n$ Toeplitz matrix $(a_{j-k})_{j,k=0}^{n-1}$.
We are interested in the behavior of the spectral norm $\n T_n\n$ as $n \to \iy$.
Notice that if the matrix $T_n$ is positive definite, then
$\n T_n\n$ is just the maximal eigenvalue of $T_n$.

\vsk
If there is a function $a \in L^{1}(\bT)$ such that $\{a_k\}_{k \in \bZ}$ is the sequence
of the Fourier coefficients of $a$, that is,
$a_k=\frac{1}{2\pi}\int_{-\pi}^{\pi} a(e^{i\tht})e^{-ik\tht}d\tht$,
we call $a$ the symbol of the sequence $\{T_n\}$ and denote $T_n$ by $T_n(a)$.
The case where $a$ is in $L^\iy(\bT)$ is easy, since then $\n T_n(a)\n \to \n a\n_\iy$ as $n \to \iy$.
Things are more complicated for symbols $a$ in $L^{1}(\bT) \setminus L^\iy(\bT)$. We here focus our attention
on so-called Fisher-Hartwig symbols with a single singularity, that is, we consider
functions $a$ of the form
\[a(t)=|t-t_0|^{-2\al}\ph_{\be, t_0}(t)b(t) \quad (t \in \bT),\]
where
$t_0 \in \bT$, $\al$ is a complex number subject to the constraint $0 < \Ret \al <1/2$,
$\be$ is a complex number satisfying $-1/2 < \Ret \be \le 1/2$, the function $\ph_{\be,t_0}$ is defined as
\[\ph_{\be,t_0}(t)=\exp(i\be\,{\rm arg}\,(-t/t_0)) \quad (t \in \bT)\]
with ${\rm arg}\,z \in (-\pi,\pi]$, and $b$ is a function in $L^\iy(\bT)$ that is continuous at
$t_0$ and does not vanish
at $t_0$. The hypothesis $0 < \Ret \al < 1/2$ ensures that $a \in L^{1}(\bT) \setminus L^\iy(\bT)$.
We should mention that if $a$ is any piecewise continuous function on $\bT$ with a single jump,
say at $t_0 \in \bT$, and $a(t_0\pm 0) \neq 0$, then $a$
can be written in the form $a=\ph_{\be,t_0}b$ with $-1/2 < \Ret \be \le 1/2$
and a continuous function $b$. Indeed,
since $\ph_{\be,t_0}(t_0-0)=e^{\pi i \be}$ and $\ph_{\be,t_0}(t_0+0)=e^{-\pi i \be}$,
it suffices to choose ${\rm Im}\,\be \in (-\iy,\iy)$ and ${\rm Re}\,\be \in (-1/2, 1/2]$
so that
\[\frac{a(t_0+0)}{a(t_0-0)}=e^{-2\pi i \be}=e^{2\pi \,{\rm Im}\,\be}\,e^{-2\pi i\,{\rm Re}\,\be}.\]
Our main result says that
\[\n T_n(a)\n \sim C_{\al,\be}\, n^{2\,\Ret \al}\, |b(t_0)| \quad \mbox{as} \quad n \to \iy\]
where $C_{\al,\be}$ is a completely identified constant depending only on $\al$ and $\be$
and where $x_n \sim y_n$ means that $x_n/y_n \to 1$. We will also establish results for symbols with more than
one Fisher-Hartwig singularity.

\vsk
For the exciting story behind Toeplitz matrices with Fisher-Hartwig symbols and their
determinants we refer to the books \cite{BoSiUni}, \cite{BoSiBu} and the
papers \cite{Bo}, \cite{Ehr}.
For general Toeplitz
matrices, the asymptotic distribution
of the singular values and the asymptotics of the extreme singular values
has been studied by many authors, and we allow us to abstain {from} giving
an ample list of references here.
These investigations are mainly directed to the collective distribution
of the singular values (Szeg\"o-Avram-Parter theorems) or to the behavior of
the extreme singular values of $T_n(a)$ for symbols $a$ in $L^\iy(\bT)$.
The asymptotics of the smallest singular value is governed by the
nature of the zeros of the symbol $a$. This implies that the rate of convergence of the
largest singular value ($=$ the norm) of $T_n(a)$ to $\n a\n _\iy$
depends on the zeros of the function $\n a\n_\iy -|a(t)|$. These results
are not applicable to Toeplitz matrices with symbols in $L^{1}(\bT) \setminus L^\iy(\bT)$
or to Toeplitz matrices ``without symbols.'' Such matrices are considered in
\cite{BoGru1}, \cite{BoGru2}, \cite{Tre1}, \cite{Tre2}, \cite{Tre3}, \cite{TyZa},
for example, but the focus of these papers is not on the problem we are interested in here.

\vsk
Under the sole assumption that $b$ be in $L^\iy(\bT)$, the method of \cite{BoGru1}
yields the estimate $\n T_n(a) \n \le C_2\,n^{2\,\Ret \al}$ with some finite constant $C_2$.
If $\al$ is real, $\ph_\be b$ is real-valued and ${\rm essinf}\,b >0$, one can also
proceed as in \cite{BoGru1} to show the existence of a positive constant $C_1$ such that
$\n T_n(a) \n \ge C_1\,n^{2\,\Ret \al}$. Such estimates were also derived in \cite{LuHu}
by different arguments. These two-sided bounds are useful in several contexts (see \cite{LuHu}, for example),
but they are clearly far away {from} the precise asymptotics
$\n T_n(a)\n \sim C_{\al,\be}\, n^{2\,\Ret \al}\, |b(t_0)|$.

\vsk
The approach of the present paper is based on an idea of Harold Widom \cite{Wi88}, \cite{Wi100},
\cite{Wi106}: we construct integral operators $K_n$ on $L^{2}(0,1)$ such that
$\n T_n(a)\n =n^{2\,\Ret\al} \,\n K_n\n$ and prove that $K_n$ converges to
some integral operator $K$ in the operator norm on on $L^{2}(0,1)$,
which implies that $\n K_n \n \to \n K\n$.

\vsk
For nonnegative symbols, the results of this paper are of interest in the analysis of time series
with long memory. The $n$th covariance matrix of a time series is a positive definite
Toeplitz matrix $T_n(a)=(a_{j-k})_{j,k=1}^n$ and one wants to know its largest eigenvalue.
If the series has a short memory, then $a_n$
goes rapidly to zero as $|n| \to \iy$ and hence $\{a_n\}$ is the sequence of the Fourier coefficients
of a function $a \in L^\iy(\bT)$. However, in the case of a long range memory,
the numbers $a_n$ may be of
the order $|n|^{2\al-1}$ ($0<\al<1/2$), which leads to symbols $a \in L^{1}(\bT)\setminus L^\iy(\bT)$.
The symbol $a(t)=|t-t_0|^{-2\,\al}\,b(t)$ is especially popular and will be considered in detail in
Section~\ref{S5}. For more on Toeplitz matrices in time series we refer the reader to
\cite{BroDa}, \cite{Dah}, \cite{Dou}, \cite{LeRe}, \cite{LuHu}.

\vsk
{\bf Acknowledgement.} We thank Fanny Godet and Bernhard Beckermann for drawing our attention
to the problem considered in this paper. We are also greatly indebted to Torsten Ehrhardt for
many valuable remarks.

\section{A Special Class of Toeplitz Matrices} \label{Selem}

We begin with a simple observation.

\begin{prop} \label{Prop El1}
Let $\ga$ be a real number.
If $|a_{\pm n}| = O(n^\ga)$ as $n \to \iy$, then $\n T_n \n$ converges to a finite limit
for $\ga < -1$, $\n T_n \n =O(\log n)$ for $\ga = -1$, $\n T_n \n =O(n^{\ga+1})$
for $\ga > -1$. If $|a_{\pm n}| = o(n^\ga)$ as $n \to \iy$, then $\n T_n \n =o(\log n)$
for $\ga = -1$ and $\n T_n \n =o(n^{\ga+1})$
for $\ga > -1$.
\end{prop}

{\em Proof.} In the case $\ga < -1$, the sequence $\{a_k\}$ is the sequence of the Fourier
coefficients of a continuous function $a$ and hence $\n T_n \n = \n T_n(a) \n \to
\n a\n_\iy$. The spectral norm of a Toeplitz matrix one diagonal of which is occupied by
units and the remaining diagonals of which are zero equals $1$. This implies that
$\n T_n\n \le \sum_{k=-(n-1)}^{n-1} |a_k|$ and therefore yields the assertions
concerning $\ga =-1$ and $\ga >-1$. $\quad \square$

\vsg
Let $A_n=(a_{j,k})_{j,k=0}^{n-1}$ be an $n \times n$ matrix with complex entries.
We denote by $G_n$ the integral operator on $L^{2}(0,1)$ with the kernel
\[g_n(x,y)=a_{[nx],[ny]}, \quad (x,y) \in (0,1)^{2},\]
where $[\xi]$ denotes the integral part of $\xi$.

\begin{lem} {\rm \sc (Widom)}\label{Lem 2.1}
The spectral norm of $A_n$ and the operator norm of $G_n$ are related by the equality
$\n A_n\n = n \n G_n\n$.
\end{lem}

{\em Proof.} Put $I_k=(k/n,(k+1)/n)$
and consider the operators
\[S_n:\{x_k\}_{k=0}^{n-1} \mapsto
\sqrt{n}\,\sum_{k=0}^{n-1}x_k\chi_{I_k}, \quad
T_n: f \mapsto \left\{\sqrt{n}\int_{I_k}f(x)dx\right\}_{k=0}^{n-1}.\]
It is easily seen that $\n S_n\n =\n T_n\n =1$ and that $T_nS_n$ is the
identity operator on $\bC^n$. Since $S_nA_nT_n=nG_n$ and thus $A_n=nT_nG_nS_n$
we obtain that $\n A_n \n \ge n \n G_n\n$ and $\n A_n\n \le n \n G_n\n$.
$\quad \square$

\vsg
Let $C^+$, $C^-$, $\ga$ be complex numbers, let ${\rm Re}\,\ga >-1$,
let
$a_{\pm n} =C^\pm n^\ga$ for $n \ge 1$, and let $a_0$ be any complex number. Denote by $K_n$
and $K$ the integral operators on $L^{2}(0,1)$ with the kernels
\[k_n(x,y)=n^{-\ga}\,a_{[nx]-[ny]}\quad
\mbox{and}\quad
k(x,y)=\left\{\begin{array}{lll}
C^+\,(x-y)^\ga & \mbox{for} & x>y,\\
C^-(y-x)^\ga & \mbox{for} & x<y,\end{array}\right.\]
respectively.

\begin{lem} \label{Lem El2}
The operators $K_n$ converge to $K$ in the operator norm on $L^{2}(0,1)$.
\end{lem}

{\em Proof.} Fix a $\mu \in (0,1)$ sufficiently close to $1$
such that $(1-\mu)|\Ret \ga| < \mu$ and $2\,\mu\,\Ret \ga < 1+2\,\Ret \ga$. Put
\begin{eqnarray*}
& & k_n^1(x,y)=\left\{\begin{array}{ll} k(x,y) & \mbox{if} \;\: |x-y|>n^{\mu-1},\\
0 & \mbox{otherwise}, \end{array}\right.
\; k_n^{2}(x,y)=\left\{\begin{array}{ll} k(x,y) & \mbox{if} \;\: |x-y|<n^{\mu-1},\\
0 & \mbox{otherwise}, \end{array}\right.\\
& & \ell_n^1(x,y)=\left\{\begin{array}{ll} k_n(x,y) & \!\!\mbox{if} \;\: |x-y|>n^{\mu-1},\\
0 & \!\!\mbox{otherwise}, \end{array}\right.
\; \ell_n^{2}(x,y)=\left\{\begin{array}{ll} k_n(x,y) & \!\!\mbox{if} \;\: |x-y|<n^{\mu-1},\\
0 & \!\!\mbox{otherwise}, \end{array}\right.
\end{eqnarray*}
and denote by $K_n^1, K_n^{2}, L_n^1, L_n^{2}$ the integral operators on $L^{2}(0,1)$ with the
kernels $k_n^1, k_n^{2}, \ell_n^1, \ell_n^{2}$, respectively. We have $K=K_n^1+K_n^{2}$ and
$K_n=L_n^1+L_n^{2}$. Thus,
\[\n K-K_n\n \le \n K_n^1-L_n^1\n +\n K_n^{2}\n +\n L_n^{2}\n.\]
We show that each term on the right goes to zero as $n \to \iy$.

\vsk
To prove that $\n K^1-K_n^1\n \to 0$ it suffices to show that $|k_n^1(x,y)-\ell_n^1(x,y)|$
converges uniformly to zero for $|x-y|>n^{\mu-1}$. We may assume that $x>y$, since the case
$x<y$ can be tackled analogously. Thus, let $x-y>n^{\mu-1}$. As $[nx]-[ny] =n(x-y)+\eps_n$
with $|\eps_n|=|\eps_n(x,y)| \le 2$, we get
\begin{eqnarray*}
\ell_n^1(x,y) & = & C^+n^{-\ga}([nx]-[ny])^\ga
=C^+n^{-\ga}(n(x-y)+\eps_n)^\ga\\
& = & C^+(x-y)^\ga \left(1+\frac{\eps_n}{n(x-y)}\right)^\ga.
\end{eqnarray*}
Since $n(x-y) > n^\mu$, it follows that
$\ell_n^1(x,y)=C^+(x-y)^\ga(1+O(n^{-\mu}))$
uniformly in $x$ and $y$. Hence
\[|k_n^1(x,y)-\ell_n^1(x,y)| =|(x-y)^\ga|\,O(n^{-\mu})=(x-y)^{\Ret \ga}\,O(n^{-\mu})\]
uniformly in $x$ and $y$. If $\Ret \ga \ge 0$, this goes to zero uniformly in $x$ and $y$.
In the case where $\Ret \ga <0$, we use the inequality $x-y > n^{\mu-1}$ to obtain that
\[(x-y)^{\Ret \ga}\,O(n^{-\mu})=O\left(n^{(1-\mu)|\Ret \ga|}n^{-\mu}\right)\]
uniformly in $x$ and $y$, which is $o(1)$ because $(1-\mu)|\Ret \ga| < \mu$. We so have
proved that $\n K^1-K_n^1\n \to 0$ as $n \to \iy$.

\vsk
The operator $K_n^{2}$ is the compression to $L^{2}(0,1)$ of the operator of convolution
on $L^{2}(\bR)$ by the kernel
\[\kappa(x)=\left\{\begin{array}{ll} C^+x^\ga & \mbox{for} \;\: 0<x<n^{\mu-1},\\
C^-|x|^\ga & \mbox{for} \;\: -n^{\mu-1}<x<0,\\
0 & \mbox{for} \;\: |x| >n^{\mu-1}.\end{array}\right.\]
The norm of a convolution operator on $L^{2}(\bR)$ is the maximum of the modulus of
the Fourier transform
\[(F\kappa)(\xi)=\int_{\bR} \kappa(x)e^{i\xi x}dx \quad (\xi \in \bR)\]
of its convolution kernel $\kappa(x)$. Hence
\begin{eqnarray*}
\n K_n^{2}\n & \le & \max_{\xi \in \bR} |(F\kappa)(\xi)| \le \int_{\bR} |\kappa(x)|dx\\
& = & \int_{-n^{\mu-1}}^1 C^-|x|^{\Ret \ga} dx+ \int_0^{n^{\mu-1}} C^+x^{\Ret \ga} dx
= O\left(n^{(\mu-1)(\Ret \ga+1)}\right),
\end{eqnarray*}
which proves that $\n K_n^{2}\n \to 0$ as $n \to \iy$.

\vsk
Let us consider the norm $\n L_n^{2}\n$. The kernel $\ell_n^{2}(x,y)$ is supported in the strip
$|x-y|<n^{\mu-1}$. Let $\tilde{\ell}_n^{2}(x,y)$ be $k_n(x,y)$ for $(x,y)$ in the
staircase-like bordered strip $|[nx]-[ny]| < n^{\mu}$ and be zero otherwise. Denote by
$\tilde{L}_n^{2}$ the corresponding integral operator. The difference $\ell_n^{2}(x,y)-
\tilde{\ell}_n^{2}(x,y)$ is supported in about $4(n-n^\mu)=O(n)$ squares of side length
$1/n$, and in these squares the absolute value of the difference is about $n^{-\ga}a_{\pm[n^\mu]}
=O(n^{-\Ret \ga}n^{\mu\Ret \ga})$.
Consequently, the squared
Hilbert-Schmidt norm $\n L_n^{2}-\tilde{L}_n^{2}\n_2^{2}$ is at most a constant times
$n\,n^{-2\,\Ret \ga}n^{2\,\mu\Ret \ga} (1/n)^{2}$,
which goes to zero because
$1-2\,\Ret \ga +2\,\mu\,\Ret\ga-2= 2\,\mu\,\Ret\ga-(1+2\,\Ret\ga)<0$.
We are therefore left with proving that $\n \tilde{L}_n^{2}\n \to 0$. Let $T_n=(b_{j-k})_{j,k=0}^{n-1}$
where $b_k=a_k$ for $|k| \le n^\mu$ and $b_k=0$ otherwise. Lemma \ref{Lem 2.1} implies that
$\n \tilde{L}_n^{2}\n =(1/n)n^{-\Ret \ga} \n T_n\n$, and since
\[\n T_n\n \le \sum_{k=-n^\mu}^{n^\mu} |b_k| =O\left(\sum_{k=-n^\mu}^{n^\mu} k^{\Ret \ga}\right)
=O\left(n^{\mu(\Ret \ga+1)}\right),\]
we finally get $\n \tilde{L}_n^{2}\n =O\left(n^{(\mu-1)(\Ret \ga+1)}\right)=o(1)$.
$\quad \square$

\begin{thm}\label{Theo El3}
Let $T_n=(a_{j-k})_{j,k=0}^{n-1}$ where $a_{\pm n} =C^\pm n^\ga (1+o(1))$ as $n \to \iy$
with complex numbers $C^+$, $C^-$, $\ga$ such that ${\rm Re}\,\ga >-1$ and at least one of the numbers
$C^+$ and $C^-$ is nonzero.
Then
\[\n T_n \n \sim \n K\n \, n^{\Ret \ga+1},\]
where $K$ is the integral operator on $L^{2}(0,1)$ whose kernel is $C^+(x-y)^\ga$ for $x>y$
and $C^-(y-x)^\ga$ for $x<y$.
\end{thm}

{\em Proof.} Write $T_n=S_n+D_n$ with $S_n=(b_{j-k})_{j,k=0}^{n-1}$,
$D_n=(d_{j-k})_{j,k=0}^{n-1}$, $b_{\pm n}=C^\pm n^\ga$, $d_{\pm n}=o(n^\ga)$.
{From} Lemma \ref{Lem 2.1} we deduce that
$\n S_n\n/n$ equals the norm of the integral operator $n^{\ga}\,K_n$
where $K_n$ has the kernel $n^{-\ga}b_{[nx]-[ny]}$. Lemma \ref{Lem El2}
implies that $\n K_n -K\n \to 0$ and thus $\n K_n\n \to \n K\n$.
Consequently, \[\n S_n\n =\n K\n n^{\Ret \ga +1}(1+o(1)).\]
Proposition \ref{Prop El1} yields $\n D_n\n=o(n^{\Ret \ga +1})$.
$\quad \square$

\begin{thm}\label{Theo El4}
Let $B_1^+, \ldots, B_Q^+, B_1^-, \ldots, B_Q^-$ be complex numbers and suppose
at least one of these numbers is nonzero. Let further $\ga_1, \ldots, \ga_Q$ be complex numbers
such that
$\Ret \ga_s >-1$ for all $s$.
Put $\om=e^{2\pi i/Q}$.
Let $T_n=(a_{j-k})_{j,k=0}^{n-1}$ where
\[a_{\pm n} =\sum_{s=1}^Q B_s^\pm \om^{\pm sn} n^{\ga_s} (1+o(1))\;\: \mbox{as}\;\:n \to \iy.\]
Put $\Ret \ga :=\max_s \Ret \ga_s$ and let $S=\{s: \Ret \ga_s=\Ret \ga\}$.Then
\[\n T_n \n \sim \max_{s\in S} \n K_s\n \, n^{\Ret \ga+1},\]
where $K_s$ is the integral operator on $L^{2}(0,1)$
whose kernel is $B_s^+(x-y)^{\ga_s}$ for $x>y$
and $B_s^-(y-x)^{\ga_s}$ for $x<y$.
\end{thm}

{\em Proof.} Assume first that $n=mQ$ with a natural number $m$. We rearrange the rows
of $T_{mQ}$ by first taking the rows $1, Q+1, 2Q+1, \ldots$, then the rows
$2, Q+2, 2Q+2, \ldots$, and so on. Then we make the same rearrangement with the columns.
The resulting matrix has the same spectral norm as $T_{mQ}$ and is a block Toeplitz matrix
$(A_{j-k})_{j,k=0}^{Q-1}$ whose blocks are the Toeplitz matrices
$A_k=(a_{k+(u-v)Q})_{u,v=0}^{m-1}$.
For $0 \le |k| \le Q-1$, let $D_k$ be the Toeplitz matrix
$D_k=(d^{(k)}_{u-v})_{u,v=0}^{m-1}$ given by $d^{(k)}_0=0$ and
\[d^{(k)}_{\pm \nu}=\sum_{s=1}^Q B_s^\pm \om^{sk}(\nu Q)^{\ga_s}\quad (\nu \ge 1).\]
If $\nu \to \infty$, then eventually $k+\nu Q \ge 1$ and hence
\[a_{k+\nu Q}-d^{(k)}_\nu =\sum_{s=1}^Q B_s^+\om^{sk}
\left[(k+\nu Q)^{\ga_s}(1+o(1))-(\nu Q)^{\ga_s}\right].\]
The modulus of the term in brackets is
\begin{eqnarray*}
& & (\nu Q)^{\Ret \ga_s}\left|\left(1+\frac{k}{\nu Q}\right)^{\ga_s}(1+o(1))-1\right|\\
& & = (\nu Q)^{\Ret \ga_s} |(1+o(1))(1+o(1))-1|=(\nu Q)^{\Ret \ga_s}\,o(1)
=o(\nu^{\Ret \ga_s}).
\end{eqnarray*}
An analogous estimate holds for $\nu \to -\iy$. {From} Proposition \ref{Prop El1}
we therefore deduce that $A_k=D_k+E_k$ with $\n E_k\n =o(m^{\Ret \ga +1})$.
It follows that
\[\n T_{mQ}\n =\n (D_{j-k})_{j,k=0}^{Q-1}\n +o(m^{\Ret \ga+1}).\]
Now, for $s \in \{1, \ldots, Q\}$, put $H_s=(h^{(s)}_{u-v})_{u,v=0}^{m-1}$ where $h_0^{(s)}=0$ and
$h^{(s)}_{\pm \nu}=B_s^\pm (\nu Q)^{\ga_s}$ for $\nu \ge 1$. Then
\[D_{j-k}=\sum_{s=1}^Q \om^{s(j-k)}H_s.\]
Let $F_\pm$ be the block Fourier matrices
$F_\pm=(\om^{\pm jk}I_{m \times m})_{j,k=1}^{Q}$.
The $j,k$ block entry of $F_+\,{\rm diag}\,(H_1, \ldots, H_Q)\,F_-$ equals
$\sum_{s=1}^Q \om^{js}H_s\om^{-sk}=D_{j-k}$.
Consequently,
\[(D_{j-k})_{j,k=0}^{Q-1}=(D_{j-k})_{j,k=1}^Q=F_+\,{\rm diag}\,(H_1, \ldots, H_Q)\,F_-.\]
Since the matrices $(1/\sqrt{Q})\,F_\pm$ are unitary, we get
\[\n (D_{j-k})_{j,k=0}^{Q-1}\n =\sqrt{Q}\,\,\n\, {\rm diag}\,(H_1, \ldots, H_Q)\,\n \,\sqrt{Q}
= Q\,\max_s \n H_s\n.\]
Theorem \ref{Theo El3} gives
\[\n H_s\n = m^{\Ret \ga_s}\,Q^{\Ret \ga_s+1}\,\n K_s\n\,(1+o(1)).\]
In summary,
\begin{eqnarray*}
\n T_{mQ}\n & = & Q\,\max_s\, m^{\Ret \ga_s+1}\,Q^{\Ret \ga_s}\,\n K_s\n\,(1+o(1))+o(m^{\Ret \ga+1})\\
& = & (m Q)^{\Ret \ga+1}\, \max_{s \in S} \n K_s\n +o(m^{\Ret \ga+1})\\
& \sim & (m Q)^{\Ret \ga+1}\, \max_{s \in S} \n K_s\n.
\end{eqnarray*}
Finally, if $n$ is not divisible by $Q$, we can obtain $T_n$ {from} $T_{mQ}$ by adding at
most $Q-1$ rows and columns. The spectral norm of a matrix with a single nonzero row or column
is the $\ell^{2}$ norm of this row or column, which in the case at hand does not exceed the
square root of
$\sum_{j=-(n-1)}^{n-1}|a_j|^{2}=O(n^{2\Ret \ga +1})$,
that is, $O(n^{\Ret \ga +1/2})=o(n^{\Ret \ga +1})$. This completes the proof. $\quad \square$

\section{A Single Fisher-Hartwig Singularity}\label{S2}

We first consider the pure Fisher-Hartwig singularity at $t_0=1$, that is, the function
\[\si(t)=|t-1|^{-2\al}\ph_{\be,1}(t)\]
with $0 < \Ret \al < 1/2$ and $-1/2 < \Ret \be \le 1/2$.
The Fourier coefficients of $\si$ are
\[\si_n=(-1)^n\,\frac{\Ga(1-2\al)}{\Ga(-\al+\be+1-n)\Ga(-\al-\be+1+n)},\]
with the convention that $\si_n:=0$ for $n <0$ if $\al=-\be$ and $\si_n:=0$ for $n>0$ if $\al=\be$
(see \cite[Lemma 6.18]{BoSiBu}). Using the formula
\[\Ga(1-z)=\frac{\pi z}{\sin \pi z}\, \frac{1}{\Ga(1+z)}\]
we see that
\begin{eqnarray*}
\si_n & = & (-1)^n \,\Ga(1-2\al)\,\frac{\sin \pi(n+\al-\be)}{\pi(n+\al-\be)}\,
\frac{\Ga(n+1+\al-\be)}{\Ga(n+1-\al-\be)}\\
& = & \Ga(1-2\al)\,\frac{\sin \pi(\al-\be)}{\pi(n+\al-\be)}\,
\frac{\Ga(n+1+\al-\be)}{\Ga(n+1-\al-\be)}
\end{eqnarray*}
for $n \ge 0$ and
\begin{eqnarray*}
\si_{-n} & = & (-1)^n\,\frac{\Ga(1-2\al)}{\Ga(-\al+\be+1+n)\Ga(-\al-\be+1-n)}\\
& = & (-1)^n \,\Ga(1-2\al)\,\frac{\sin \pi(n+\al+\be)}{\pi(n+\al+\be)}\,
\frac{\Ga(n+1+\al+\be)}{\Ga(n+1-\al+\be)}\\
& = & \Ga(1-2\al)\,\frac{\sin \pi(\al+\be)}{\pi(n+\al+\be)}\,
\frac{\Ga(n+1+\al+\be)}{\Ga(n+1-\al+\be)}
\end{eqnarray*}
for $n \ge 0$. The asymptotic formula
$\Ga(n+\ga)/\Ga(n+\de)\sim n^{\ga-\de}$
($n \to \iy$)
shows that
\[\si_n=C_{\al,\be}^+\,n^{2\al-1} (1+o(1)),
\quad \si_{-n}=C_{\al,\be}^-\,n^{2\al-1}(1+o(1))\]
as $n \to \iy$, where
\[C_{\al,\be}^\pm = \Ga(1-2\al)\,\frac{\sin \pi(\al \mp \be)}{\pi}.\]
We denote by $K$ the integral operator on $L^{2}(0,1)$ with the kernel
\[k(x,y)=\left\{\begin{array}{lll}
C_{\al,\be}^+\,(x-y)^{2\al-1} & \mbox{for} & x>y,\\
C_{\al,\be}^-\,(y-x)^{2\al-1} & \mbox{for} & x<y.\end{array}\right.\]
Obviously, $\n K \n >0$.

\begin{thm} \label{Theo 2.4}
Suppose $\si(t)=|t-t_0|^{-2\al}\ph_{\be,t_0}(t)$ with $t_0 \in \bT$, $0 < \Ret \al < 1/2$,
$-1/2 < \Ret \be \le 1/2$. Then
\[\n T_n(\si)\n \sim \n K \n \,n^{2\,\Ret\al}.\]
\end{thm}

{\em Proof.} We slightly change notation and denote the function $\si$ defined as
$\si(t)=|t-1|^{-2\al}\ph_{\be,1}(t)$ by $\si^0$. The $\si$ of the present theorem
results {from} $\si^0$ by replacing $t_0=1$ with a general
$t_0 \in \bT$. The only change in the Fourier coefficients is that the $(-1)^n$ in $(\si^0)_n$
becomes
$(-1/t_0)^n$ in $\si_n$ and hence $T_n(\si)=\Lambda\, T_n(\si^0)\,\Lambda^{-1}$
where $\Lambda:={\rm diag}\,(1,t_0^{-1}, \ldots,t_0^{-(n-1)})$. Therefore
$\n T_n(\si)\n =\n T_n(\si^0)\n$. Taking into account that $(\si^0)_{\pm n}
=C^\pm_{\al,\be}\, n^{2\,\al-1}\,(1+o(1))$ and using Theorem \ref{Theo El3}
we arrive at the desired formula. $\quad \square$

\begin{prop} \label{Prop C}
If $\si$ is as in Theorem \ref{Theo 2.4} and $c \in L^\iy(\bT)$ is continuous and zero at $t_0$,
then
\[\n T_n(\si c)\n = o(n^{2\,\Ret \al}).\]
\end{prop}

{\em Proof.} Without loss of generality assume that $t_0=1$. Writing $\si = \Ret \si+i\,{\rm Im}\,\si$
and $c=\Ret c +i \,{\rm Im}\,c$ we get $T_n(\si c)=T_n(\Ret \si\,\Ret c)+ \ldots$ (four terms) and thus
$\n T_n(\si c)\n \le \n T_n(\Ret \si\,\Ret c)\n + \ldots$. The matrix $T_n(\Ret \si\,\Ret c)$
is Hermitian and hence
\begin{eqnarray*}
\n T_n(\Ret \si\,\Ret c)\n & = & \max_{\psi \in \bC^n \setminus\{0\}}
\frac{|(T_n(\Ret \si\,\Ret c)\psi,\psi)|}{\n \psi \n^{2}}\\
& = & \max_{\ph \in {\mathcal P}_n\setminus\{0\}} \frac{1}{\n \ph \n^{2}}\,
\left|\,\frac{1}{2\pi}\int_{-\pi}^\pi \Ret \si(x) \,\Ret c(x)\,|\ph(x)|^{2}dx\right|\\
& \le & \max_{\ph \in {\mathcal P}_n\setminus\{0\}} \frac{1}{\n \ph \n^{2}}\,
\frac{1}{2\pi}\int_{-\pi}^\pi |\Ret \si(x)| \,|\Ret c(x)|\,|\ph(x)|^{2}dx,
\end{eqnarray*}
where ${\mathcal P}_n$ is the set of all trigonometric polynomials of the form
$\ph(x)=\ph_0+\ph_1e^{ix}+\ldots+\ph_{n-1}e^{i(n-1)x}$. Notice that
\begin{eqnarray*}
\n \ph\n^{2}_\iy & = & \n \ph_0+\ph_1e^{ix}+\ldots+\ph_{n-1}e^{i(n-1)x}\n_\iy^{2}
\le (|\ph_0|+|\ph_1|+\ldots+|\ph_{n-1}|)^{2}\\
& \le & n\,(|\ph_0|^{2}+|\ph_1|^{2}+\ldots+|\ph_{n-1}|^{2})=\frac{n}{2\pi}\int_{-\pi}^\pi
|\ph(x)|^{2}dx=\frac{n}{2\pi}\,\n\ph\n^{2}.
\end{eqnarray*}
Clearly,
\[\Ret \si(x)=\left|2\,\sin\frac{x}{2}\right|^{-2\,\Ret \al}
\cos\left({\rm Im}\,\al \,\log\left|2\sin\frac{x}{2}\right|\,\right)\,|\ph_{\be,1}(x)|,\]
which is $O(|x|^{-2\,\Ret \al})$ as $x \to 0$.
We split the integral into $\int_{|x| < \pi/n}$, $\int_{\pi/n < |x|<\pi/\sqrt{n}}\,$,
and $\int_{\pi/\sqrt{n}<|x|<\pi}$. The integral over $|x|<\pi/n$ is at most
\begin{eqnarray*}
& & C_1\, \sup_{|x|<\pi/n} |\Ret c(x)|\,\n \ph\n_\iy^{2}\, \int_{|x|<\pi/n} |x|^{-2\,\Ret \al}dx\\
& & \le C_2\, \sup_{|x|<\pi/n} |\Ret c(x)|\,\frac{n}{2\pi}\,\n \ph\n^{2}\,n^{2\,\Ret \al-1}
=o(n^{2\,\Ret \al})\, \n \ph\n^{2}
\end{eqnarray*}
because $\Ret c(x) \to 0$ as $x \to 0$; here $\sup$ means ${\rm esssup}$. The integral over
the interval $\pi/n < |x|<\pi/\sqrt{n}$ has the upper bound
\begin{eqnarray*}
& & C_3\,\sup_{|x|<\pi/\sqrt{n}}|\Ret c(x)| \,\int_{|x|>\pi/n}
|x|^{-2\,\Ret \al}|\ph(x)|^{2}dx\\
& & \le C_3\,\sup_{|x|<\pi/\sqrt{n}}|\Ret c(x)|\,\frac{n^{2\,\Ret \al}}{\pi^{2\,\Ret \al}}\,
\int_{|x|>\pi/n}|\ph(x)|^{2}dx\\
& & \le C_3\,\sup_{|x|<\pi/\sqrt{n}}|\Ret c(x)|\,\frac{n^{2\,\Ret \al}}{\pi^{2\,\Ret \al}}\,
\n \ph\n^{2} =o(n^{2\Ret \al})\,\n \ph\n^{2},
\end{eqnarray*}
again because $\Ret c(x) \to 0$ as $x \to 0$. Finally, the integral over $|x| >\pi/\sqrt{n}$
does not exceed
\begin{eqnarray*}
& & C_4\,\n \Ret c\n_\iy \,\int_{|x|>\pi/\sqrt{n}} |x|^{-2\,\Ret \al}|\ph(x)|^{2}dx\\
& & \le C_4\,\n \Ret c\n_\iy\,\frac{n^{\Ret \al}}{\pi^{2\,\Ret \al}}
\int_{|x|>\pi/\sqrt{n}}|\ph(x)|^{2}dx\\
& & \le C_4\,\n \Ret c\n_\iy\,\frac{n^{\Ret \al}}{\pi^{2\,\Ret \al}}\,\n\ph\n^{2}
=O(n^{\Ret \al})\,\n \ph \n^{2}=o(n^{2\,\Ret \al}) \,\n \ph \n^{2}.
\end{eqnarray*}
This proves that $\n T_n(\Ret \si\,\Ret c)\n = o(n^{2\,\Ret \al})$. Analogously one can show that
\[\n T_n(\Ret \si\,{\rm Im}\, c)\n, \quad \n T_n({\rm Im}\, \si\,{\rm Re}\, c)\n, \quad
\n T_n({\rm Im}\, \si\,{\rm Im}\, c)\n\]
are $o(n^{2\,\Ret \al})$. $\quad \square$

\begin{thm} \label{Theo 2.5}
Let $a=\si b$ where $\si$ is as in Theorem \ref{Theo 2.4} and $b$ is a function in $L^\iy(\bT)$
that is continuous at $t_0$ and does not vanish at $t_0$. Then
\[\n T_n(a)\n \sim \n K\n\,|b(t_0)|\, n^{2\,\Ret \al}.\]
\end{thm}

{\em Proof.}  We have $b(t)=b(t_0)+c(t)$ with $b(t_0) \neq 0$ and a function $c \in L^\iy(\bT)$
that is continuous and zero at $t_0$. It follows that
$T_n(a)=b(t_0)\,T_n(\si)+T_n(\si c)$.
Theorem~\ref{Theo 2.4} yields
\[\n b(t_0)\,T_n(\si)\n=|b(t_0)|\,\n T_n(\si)\n =|b(t_0)|\, \n K\n\,n^{2\,\Ret \al}\,(1+o(1)),\]
and Proposition \ref{Prop C} gives $\n T_n(\si c)\n =o(n^{2\,\Ret \al})$.
$\quad \square$

\section{Several Fisher-Hartwig Singularities} \label{S3}

Let $R \ge 2$ and
\[a(t)=b(t)\prod_{r=1}^R |t-t_r|^{-2\al_r}\ph_{\be_r,t_r}(t) \quad (t \in \bT)\]
where $t_1, \ldots, t_R$ are distinct points on $\bT$,
$0 < \Ret \al_r< 1/2$, $-1/2 < \Ret \be_r \le 1/2$, $b \in L^\iy(\bT)$,
$b$ is continuous at the points $t_1, \ldots, t_R$, and
$b(t_r) \neq 0$ for all $r$. It is easily seen that $a$ can be written in the form
\[a(t)=\sum_{r=1}^R |t-t_r|^{-2\al_r}\ph_{\be_r,t_r}(t)\, b_r(t) \quad (t \in \bT)\]
with functions $b_r \in L^\iy(\bT)$ such that $b_r$ is continuous at $t_r$ and satisfies
$b_r(t_r) \neq 0$. Let
\[\Ret \al :=\max\{\Ret\al_1, \ldots, \Ret\al_R\}, \quad M=\{ r: \Ret \al_r =\Ret \al\}.\]
If there is only one $r_0$ such that $\Ret \al_{r_0}=\Ret \al$, then Theorem \ref{Theo 2.5}
implies that
\[\n T_n(a) \n \sim \n K_{r_0} \n \, |b_{r_0}(t_{r_0})|\,n^{2\,\Ret \al},\]
where $K_{r}$ denotes the integral operator on $L^{2}(0,1)$ associated with
$|t-t_r|^{-2\al_r}\ph_{\be_r,t_r}(t)$, that is, the integral operator whose kernel
is $C_{\al_r,\be_r}^+(x-y)^{2\,\al_r-1}$ for $x>y$ and equals
$C_{\al_r,\be_r}^-(y-x)^{2\,\al_r-1}$ for $x<y$.
The case where the maximum is attained at more than one $r$ is more difficult.

\begin{con} \label{Con 4.1}
We have
\[\n T_n(a)\n \sim \max_{r \in M} \n K_r\n \,|b(t_r)|\, n^{2\,\Ret \al}.\]
\end{con}

The following result confirms this conjecture in a sufficiently interesting special case.

\begin{thm} \label{Theo 4.2}
If there is a $t_0 \in \bT$ such that, for every $r$, $t_r=e^{2\pi i \ph_r}t_0$ with a rational number $\ph_r$,
then Conjecture \ref{Con 4.1} is true.
\end{thm}

{\em Proof.} As passage {from} $a(t)$ to $a(t/t_0)$ does not change the spectral norm of the Toeplitz matrix
(recall the proof of Theorem \ref{Theo 2.4}), we may without loss of generality assume that $t_0=1$.
Put $\si_{\al,\be,\tau}(t)=|t-\tau|^{-2\al}\ph_{\be,\tau}(t)$. The Fourier coefficients of $a$ are
\[a_n=\sum_{r=1}^R (\si_{\al_r, \be_r,t_r}\,b_r)_n= \sum_{r=1}^R
b_r(t_r)\,(\si_{\al_r, \be_r,t_r})_n +f_n\]
where $\{f_n\}$ is the sequence of the Fourier coefficients of a function $f \in L^{1}(\bT)$
for which $\n T_n(f)\n =o(n^{2\,\Ret \al})$ (Proposition~\ref{Prop C}). Furthermore,
$(\si_{\al_r, \be_r,t_r})_n =t_r^{-n}(\si_{\al_r, \be_r,1})_n$ (see once more the proof of
Theorem \ref{Theo 2.4}). Thus,
\[a_n=\sum_{r=1}^R b_r(t_r)\,t_r^{-n}\,(\si_{\al_r, \be_r,1})_n +f_n.\]
Let $t_r^{-1}=e^{2\pi ip_r/q_r}$ with a rational number $p_r/q_r \in (0,1]$ and denote by $Q$ the
least common multiple of $q_1, \ldots, q_R$. Put $\om =e^{2\pi i/Q}$. Then each $t_r^{-1}$
is of the form $\om^{k_r}$ with some $k_r\in \{1,2, \ldots, Q\}$. It follows that
\[a_n=\sum_{r=1}^R b_r(t_r) \,\om^{k_rn}\,(\si_{\al_r, \be_r,1})_n +f_n\]
with different $k_1, \ldots, k_R$ belonging to $\{1, 2, \ldots, Q\}$. {From} Section \ref{S2}
we know that
\[(\si_{\al, \be,1})_{\pm n}=C_{\al,\be}^\pm\,n^{2\al-1}(1+o(1)).\]
Hence
\[a_{\pm n}=\sum_{r=1}^R b_r(t_r)\,\om^{\pm k_rn}\,C_{\al_r,\be_r}^\pm\,n^{2\al_r-1}(1+o(1))+f_n,\]
which can be written as
\[a_{\pm n}=\sum_{s=1}^Q B_s^\pm\, \om^{\pm sn}\,n^{\ga_s}(1+o(1))+f_n\]
with $B_{k_r}^\pm=b_r(t_r)\,C_{\al_r,\be_r}^\pm$, $\ga_{k_r}=2\al_r-1$ and $B_s^\pm=0$, $\ga_s=0$
otherwise. Theorem \ref{Theo El4} shows that the spectral norm of the Toeplitz matrix
$T_n^0$ generated by
\[a_{\pm n}^0:=\sum_{s=1}^Q B_s^\pm\, \om^{\pm sn}\,n^{\ga_s}(1+o(1))\]
satisfies
\[\n T_n^0\n \sim \max_{s \in S} \n K_s^0\n\,n^{2\, \Ret \al}
\quad \mbox{with} \quad \Ret \al :=\max_s \Ret\al_s, \quad S=\{ s: \Ret \al_s =\Ret \al\},\]
where $K_s^0$ is the operator whose kernel is $B_s^+(x-y)^{\ga_s}$ for
$x>y$ and $B_s^-(y-x)^{\ga_s}$ for $x<y$. This is equivalent to saying that
\[\n T_n^0\n \sim \max_{r \in M} |b_r(t_r)|\,\n K_r\n \,n^{2\,\Ret \al}\]
where the kernel of $K_r$ is
$C_{\al_r,\be_r}^+(x-y)^{2\,\al_r-1}$ for $x>y$ and
$C_{\al_r,\be_r}^-(y-x)^{2\,\al_r-1}$ for $x<y$.
Since $\n T_n(f)\n =o(n^{2\,\Ret \al})$, we obtain that
$\n T_n \n \sim \n T_n^0\n$. $\quad \square$

\section{A Particular Singularity} \label{S5}

We finally embark on the case where
\[a(t)=|t-t_0|^{-2\al}b(t) \quad (t \in \bT)\]
with a real number $\al \in (0,1/2)$ and a function $b \in L^\iy(\bT)$ that is
continuous and nonzero at $t_0$. Theorem \ref{Theo 2.5} gives
\[\n T_n(a)\n \sim \Ga(1-2\al)\,\frac{\sin \pi \al}{\pi}\, \n K_\al\n
\,|b(t_0)|\, n^{2\al}\]
where the kernel of $K_\al$ is $|x-y|^{2\al-1}$.

\begin{prop} \label{Prop 5.1}
We have
\[\frac{1}{2\al}\,\left(\frac{2}{4\al+1}+2\,\frac{\Ga(2\al+1)\Ga(2\al+1)}{\Ga(4\al+2)}
\right)^{1/2} \le \n K_\al \n \le \frac{1}{\al}.\]
\end{prop}

{\em Proof.} We may think of $K_\al$ as the compression to $L^{2}(0,1)$ of the convolution
operator on $L^{2}(\bR)$ whose convolution kernel $\kappa(x)$ is $|x|^{2\al-1}$ for $|x| <1$
and $0$ for $|x|>1$. As in the proof of Lemma \ref{Lem El2} we therefore see that
\[\n K_\al\n \le \max_{\xi \in \bR} |(F\kappa)(\xi)| \le \int_{\bR} |\kappa(x)|\,dx
= \int_{-1}^{1} |x|^{2\al-1}\,dx
=\frac{1}{\al}.\]
Let ${\bf 1}$ be the function which is identically $1$ on $(0,1)$.
Taking into account that  \[\n K_\al\n^{2} \ge \n K_\al {\bf 1}\n^{2}/\n {\bf 1}\n^{2}
=\n K_\al{\bf 1}\n^{2}\quad \mbox{and} \quad
(K_\al {\bf 1})(x)=\frac{1}{2\al}\,(x^{2\al}+(1-x)^{2\al}),\]
we obtain that $\n K_\al\n^{2}$ is greater than or equal to
\[\frac{1}{4\al^{2}}\,\int_0^{1}(x^{2\al}+(1-x)^{2\al})^{2}\,dx
= \frac{1}{4\al^{2}}\,\left(\frac{2}{4\al+1}+2\,\frac{\Ga(2\al+1)\Ga(2\al+1)}{\Ga(4\al+2)}
\right).\]
This proves the lower bound for $\n K_\al \n$. $\quad \square$

\begin{cor} \label{Cor 5.2}
We have $\n K_\al \n \sim 1/\al$ as $\al \to 0$ and $\n K_\al \n \sim 1$ as $\al \to 1/2$.
\end{cor}

{\em Proof.} By Proposition \ref{Prop 5.1}, $\al^{2}\,\n K_\al\n^{2} \le 1$ and
\[\liminf_{\al \to 0} \al^{2}\,\n K_\al\n^{2} \ge \frac{1}{4}\left(2+2\,\frac{\Ga(1)\Ga(1)}{\Ga(2)}\right)=1,\]
which implies that $\al\,\n K_\al\n \to 1$ as $\al \to 0$.
Thinking of $K_\al-K_{1/2}$ as the convolution operator with the convolution kernel
$|x|^{2\al-1}-1$ for $|x|<1$ and $0$ for $|x|>1$, we get
\[\n K_\al-K_{1/2}\n \le \int_{-1}^{1}(|x|^{2\al-1}-1)\,dx
=\frac{1}{\al}-2 = o(1) \;\:\mbox{as}\;\: \al \to \frac{1}{2}.\]
Thus, $\n K_\al \n \to \n K_{1/2}\n$ as $\al \to 1/2$. Since $(K_{1/2}f)(x)=\int_0^{1}f(y)\,dy$,
it is easily seen that $\n K_{1/2}\n =1$. $\quad \square$

\begin{cor} \label{Cor 5.3}
We have
\begin{eqnarray*}
& & \Ga(1-2\al)\,\frac{\sin \pi \al}{\pi}\, \n K_\al\n \sim 1\;\:\mbox{as}\;\: \al \to 0,\\
& & \Ga(1-2\al)\,\frac{\sin \pi \al}{\pi}\, \n K_\al\n \sim \frac{1}{2\pi(1/2-\al)}
\;\:\mbox{as}\;\: \al \to \frac{1}{2}.
\end{eqnarray*}
\end{cor}

{\em Proof.} The asymptotics for $\al \to 0$ is immediate {from} Corollary \ref{Cor 5.2}.
For $\al \to 1/2$, Corollary \ref{Cor 5.2} and the formulas
\[\Ga(1-2\al)\,\frac{\sin \pi \al}{\pi} \sim \frac{\Ga(1-2\al)}{\pi}
=\frac{1}{\sin 2\pi \al}\,
\frac{1}{\Ga(2\al)}\sim \frac{1}{2\pi(1/2-\al)},\]
yield the asserted asymptotics. $\quad \square$

\end{document}